\documentclass[12pt,reqno]{amsart}
\usepackage[left=80pt,right=80pt]{geometry}
\usepackage{amsmath}
\usepackage{amssymb}
\usepackage{amsthm}
\usepackage{float}

\usepackage{hyperref,url}
\hypersetup{
	colorlinks=true,
  linkcolor=black,          
  citecolor=black,         
  filecolor=black,      
  urlcolor=black           
}

\theoremstyle{definition}

\allowdisplaybreaks
\makeatletter
\newcommand{\vast}{\bBigg@{4}}
\newcommand{\Vast}{\bBigg@{5}}
\makeatother

\begin{document}

\title[Combinatorial proofs of totals of Catalan statistics]{Combinatorial proofs of totals of some statistics on Catalan words}

\author[MARK SHATTUCK]{Mark Shattuck$^{*}$}
\thanks{$^{*}$Department of Mathematics, University of Tennessee, 37996 Knoxville, TN, USA
\\
\href{mailto:mshattuc@utk.edu}{\tt mshattuc@utk.edu}}

\begin{abstract}
A Catalan word is one on the alphabet of positive integers starting with $1$ in which each subsequent letter is at most one more than its predecessor.  Let $\mathcal{C}_n$ denote the set of Catalan words of length $n$. In this paper, we give combinatorial proofs of explicit formulas for the sums of several parameter values taken over all the members of $\mathcal{C}_n$.  In particular, we find such proofs for the parameters tracking the number of symmetric or $\ell$-valleys, which was previously requested by Baril et al.  Further, we find a combinatorial explanation of a related Catalan number identity whose  proof was also requested.  To carry out our arguments, we consider corresponding statistics on Dyck paths and find the cardinality of certain sets of marked Dyck paths wherein one or more of the steps is distinguished from all others.
\bigskip

\smallskip


\end{abstract}

\maketitle

\section{Introduction}

A \emph{Catalan word} $w=w_1\cdots w_n$ is a positive integral sequence with $w_1=1$ and satisfying $w_{i+1} \leq w_i+1$ for $1 \leq i \leq n-1$.  Let $\mathcal{C}_n$ denote the set of all Catalan words of length $n$ (i.e., that contain $n$ letters).  It is well known, e.g., \cite[Exercise 80]{Stan}, that the cardinality of $\mathcal{C}_n$ is given by the $n$-th Catalan number $C_n=\frac{1}{n+1}\binom{2n}{n}$ for all $n \geq0$.  For example, if $n=4$, then we have
\begin{align*}
\mathcal{C}_4=\{1111,  1112,  1121,  1122,   1123,  1211,  1212,  1221, 1222,  1223, 1231,  1232, 1233, 1234\}
\end{align*}
and $|\mathcal{C}_4|=C_4=14$.
Catalan words have connections to several other discrete structures enumerated by the Catalan numbers, such as rooted binary trees \cite{Mak} and Dyck paths \cite{Shat}, and have arisen in the context of the exhaustive generation of Gray codes for growth-restricted sequences \cite{MVaj}.  Recently, the study of Catalan words satisfying various restrictions has been an object of ongoing research interest.  In addition to the classical avoidance problem on Catalan words (e.g., \cite{BKV1,BKV2}), the distributions of several combinatorial parameters have been studied on $\mathcal{C}_n$, including those that track the number of occurrences of various kinds of strings \cite{BHH,Shat} or that concern geometric aspects of the corresponding bargraphs \cite{Ble,CMR,MRT,MSh}.

Baril et al.\ \cite{BHH} found simple closed-form expressions for the total number of symmetric and $\ell$-valleys in all the members of $\mathcal{C}_n$ by use of generating functions and raised the question of finding combinatorial proofs. Here, we provide such proofs for the totals of these parameters and several others that were established in \cite{BHH,CMR}, where only an algebraic proof was given; see Table \ref{tab1} below.  We remark that the algebraic derivation of the formula for the total in each case entailed first finding an expression for the generating function of the corresponding distribution, differentiating this expression with respect to $q$, substituting $q=1$ and finally extracting the coefficient of $x^n$ for all $n \geq 1$ (where $q$ marks the parameter on $\mathcal{C}_n$ in the distribution).

\begin{table}[htp]
{\small\begin{center}
        \begin{tabular}{|l|l|l|}\hline
             Statistic & Total on $\mathcal{C}_n$ & Reference/OEIS \#\\\hline\hline
            &&\\[-8pt]
             symmetric valley  & $(3n-2)C_{n-1}-\frac{1}{2}\sum_{k=1}^n\binom{2k}{k}$ &\cite[Cor.\ 4.7]{BHH}\\\hline
             &&\\[-8pt]
             $\ell$-valley & $\binom{2n-2\ell-1}{n-\ell-3}$ & \cite[Cor.\ 4.2]{BHH}  \\\hline
            &&\\[-8pt]
         symmetric peak      & $\sum_{k=0}^{n-3}\binom{2k+2}{k}$ & \cite[Cor.\ 5.6]{BHH}/A057552\\\hline
            &&\\[-8pt]
            $\ell$-peak  & $\binom{2n-2\ell-1}{n-\ell-2}$ &\cite[Cor.\ 5.2]{BHH}\\\hline
            &&\\[-8pt]
             runs of descents      & $\binom{2n}{n}-\binom{2n-2}{n-1}$ & \cite[Cor.\ 3.7]{BHH}/A051924\\\hline
            &&\\[-8pt]
          runs of weak ascents      & $\binom{2n-2}{n-1}$ & \cite[Cor.\ 3.5]{BHH}/A000984\\\hline
            &&\\[-8pt]
            corner of type $hu$      & $\binom{2n-1}{n-2}$ & \cite[Cor.\ 19]{CMR}/A002054\\\hline
            &&\\[-8pt]
             corner of type $dh$      & $\binom{2n-2}{n-3}$ & \cite[Cor.\ 17]{CMR}/A002694\\\hline
            &&\\[-8pt]
            semi-perimeter      & $\frac{1}{2}\left(\binom{2n+2}{n+1}-\binom{2n}{n}\right)$ & \cite[Cor.\ 8]{CMR}/A097613\\\hline
            &&\\[-8pt]
            area      & $\frac{1}{2}\left(4^n-\binom{2n}{n}\right)$ & \cite[Cor.\ 12]{CMR}/A000346\\\hline
        \end{tabular}
    \end{center}}
\caption{Formulas for the totals of Catalan word statistics.}\label{tab1}
\end{table}

A direct argument of the formulas in Table \ref{tab1} may be realized by demonstrating that they enumerate, equivalently, occurrences of certain strings of steps in the corresponding lattice paths.  Our arguments avoid the use of generating functions and involve defining bijections between various sets of lattice paths.  Note that several of the totals on $\mathcal{C}_n$ for the statistics in Table \ref{tab1} occur as entries in the OEIS \cite{Sloane} and in fact yield apparently new combinatorial interpretations for these sequences.  We also provide a counting argument for the Catalan number identity
\begin{equation}\label{Cniden}
C_n=\sum_{k=1}^{\lfloor\frac{n+1}{2}\rfloor}\frac{1}{n}\binom{n}{k}\binom{n-k}{k-1}2^{n-2k+1}, \qquad n \geq 1,
\end{equation}
answering a further question raised by Baril et al.\ \cite{BHH}.  This identity had arisen in conjunction with an analysis of the Catalan word statistic tracking the number of runs of weak ascents, and we prove \eqref{Cniden} by considering the distribution of an equivalent lattice path statistic.

Let $\mathcal{D}_n$ denote the set of lattice paths from $(0,0)$ to $(2n,0)$ using up $u=(1,1)$ and down $d=(1,-1)$ steps that never go below the $x$-axis.  Members of $\mathcal{D}_n$ are known as \emph{Dyck} paths.  There is a simple bijection (e.g., \cite{CMR}) between the members of $\mathcal{C}_n$ and $\mathcal{D}_n$ defined as follows.  Given $\pi=\pi_1\cdots \pi_n \in \mathcal{C}_n$, let $\imath(\pi)$ denote the member of $\mathcal{D}_n$ whose $i$-th up step from the left has final height $\pi_i$ for $1 \leq i \leq n$.  Note that this uniquely determines the positions of the down steps, and hence it is seen that $\imath$ defines a bijection between $\mathcal{C}_n$ and $\mathcal{D}_n$.  For example, if $n=6$ and $\pi=123321 \in \mathcal{C}_6$, then $\imath(\pi)=u^3dud^2ud^2ud \in \mathcal{D}_6$. For several of the statistic totals on $\mathcal{C}_n$ given in Table \ref{tab1}, it will be useful to consider the total of the corresponding parameter on $\mathcal{D}_n$ under $\imath$.  More precisely, we have that these statistics on $\mathcal{C}_n$ correspond to parameters on lattice paths which track the number of occurrences of a certain string of steps (particular to the statistic on $\mathcal{C}_n$ in question), where sometimes there are restrictions as to the location of the string.

This prompts the following definition.  Given a particular sequence $\tau$ in $\{u,d\}$, let $\mathcal{D}_n(\tau)$ denote the set of marked members of $\mathcal{D}_n$ wherein a single occurrence of $\tau$ is marked.  Note that $|\mathcal{D}_n(\tau)|$ gives the total number of occurrences of $\tau$ within all of the members of $\mathcal{D}_n$.  We will also need to consider the set of (unrestricted) lattice paths to particular points in the plane.  Given a point $(a,b)$ with integral coordinates, where $a$ and $b$ are of  the same parity and $a\geq|b|$, let $\mathcal{P}_{(a,b)}$ denote the set of all lattice paths from $(0,0)$ to $(a,b)$ using $u$ and $d$ steps.  Note that $|\mathcal{P}_{(a,b)}|=\binom{a}{(a-b)/2}$ for all such $a$ and $b$. At times, it will be useful to identify the members of $\mathcal{D}_n(\tau)$ with paths in $\mathcal{P}_{(a,b)}$ for some $a$ and $b$.  See, e.g., \cite{Deutsch,Sap,Sun} and references contained therein for results concerning the distributions of strings and other parameters on $\mathcal{D}_n$.

The organization of this paper is as follows.  In the second section, we give combinatorial proofs of the formulas for the totals of the valley, peak and run statistics in Table \ref{tab1}.  As a consequence of our arguments, we obtain combinatorial explanations of the formulas stated for the totals of the $hu$ and $dh$ corner statistics.  In the third section, we provide a combinatorial interpretation of the Catalan identity \eqref{Cniden}, where we make use of an enumerative result known as \emph{Raney's} lemma.  In the final section, we find direct arguments for the formulas found in \cite{CMR} for the totals of all semi-perimeter and area values over $\mathcal{C}_n$.  In the case of the semi-perimeter, to prove the result, we show that it is closely connected to the number of ascents, and hence to the number of occurrences of $uu$ in the corresponding Dyck path under $\imath$.  For the area parameter, we make use of a certain decomposition of Dyck paths based on the relative positions of a pair of marked steps in defining a bijection with a class of lattice paths enumerated by the total area.

We will make use of the following further definitions throughout.  The \emph{height} $i$ of a step  will refer to its final height (i.e., the $y$-coordinate of its endpoint), and hence $i \geq 1$ for a $u$ and $i \geq0$ for a $d$ step within a Dyck path. Given a sequence $\alpha$ of steps in $\{u,d\}$, let $r(\alpha)$ denote the lattice path obtained from $\alpha$ by reversing the sequence of steps and replacing each $u$ with $d$ and $d$ with $u$.  For example, if $\alpha=u^2du^2d^4u$, then $r(\alpha)=du^4d^2ud^2$.  Note that the mapping $r$ defines an involution on $\mathcal{D}_n$.  Finally, a \emph{unit} within a member of $\mathcal{D}_n$ refers to a sequence of steps of the form $u\beta d$, where the $u$ starts from the $x$-axis and $\beta$ is a possibly empty Dyck path.

\section{Valley, peak and run statistic totals}

In this section, we provide combinatorial explanations of the formulas found for the valley, peak and run statistics in Table \ref{tab1}.

\subsection{Valleys}In this subsection, we treat the symmetric valley and $\ell$-valley parameters on $\mathcal{C}_n$.\medskip

\noindent\emph{Symmetric valleys:} Recall that a \emph{symmetric valley} within a member of $\mathcal{C}_n$ is a string of letters of the form $a(a-1)^\ell a$ for some $a>1$ and $\ell \geq 1$.  It is seen that a symmetric valley in $\pi \in \mathcal{C}_n$ corresponds to a sequence of steps in $\imath(\pi)$ of the form $ud(du)^\ell u$.    Note that a member of $\mathcal{D}_n(ud(du)^\ell u)$ may be obtained by inserting $d(du)^\ell u$ just after
any $u$ step within a member of $\mathcal{D}_{n-\ell-1}$ and marking the resulting string of length $2\ell+3$, provided the $u$ step was not of height one.  Further, we must have $\ell \leq n-3$, since the marked $ud(du)^\ell u$ cannot include the first or the penultimate step within a member of $\mathcal{D}_n$.  Let $\mathcal{E}_m$ denote the subset of $\mathcal{D}_m(u)$ in whose members the marked $u$ is not of height one.  Then we have that the total number of symmetric valleys in $\mathcal{C}_n$ is given by the cardinality of $\cup_{\ell=1}^{n-3}\mathcal{E}_{n-\ell-1}$.

Note that $|\mathcal{E}_m|=mC_m-(C_{m+1}-C_m)$, by subtraction and the fact that there are $C_{m+1}-C_m$ units altogether in $\mathcal{D}_m$, and hence the same number of up steps of height one.  For a quick bijective proof of this fact, let $\pi=\alpha\pi'\beta \in \mathcal{D}_m$, where $\pi'$ is a unit distinguished from all others and $\alpha,\beta$ are possibly empty Dyck paths.  Then $\pi \mapsto u\alpha d\pi'\beta$ defines a bijection with members of $\mathcal{D}_{m+1}$ containing at least two units, of which there are $C_{m+1}-C_m$, by subtraction.  We thus have
$$|\cup_{\ell=1}^{n-3}\mathcal{E}_{n-\ell-1}|=\sum_{i=2}^{n-2}[iC_i-(C_{i+1}-C_i)]=\sum_{i=2}^{n-2}iC_i-C_{n-1}+2.$$
So to establish the formula stated in Table \ref{tab1} for the total number of symmetric valleys in $\mathcal{C}_n$, we must show
$$(3n-2)C_{n-1}-\frac{1}{2}\sum_{k=1}^{n}\binom{2k}{k}=\sum_{i=2}^{n-2}iC_i-C_{n-1}+2.$$
This equality may be rewritten as
\begin{equation}\label{symvae1}
(3n-1)C_{n-1}=1+\binom{2n-1}{n}+\binom{2n-3}{n-1}+\sum_{i=1}^{n-2}\left[\binom{2i}{i-1}+\binom{2i-1}{i}\right], \quad n \geq 3,
\end{equation}
by the facts $\frac{1}{2}\binom{2i}{i}=\binom{2i-1}{i}$ and $iC_i=\binom{2i}{i-1}$, both of which are readily explained bijectively.

Upon marking one of the points along the path (including possibly the starting point) or any up step or the final down step within a member of $\mathcal{C}_{n-1}$, one can show
$$(3n-1)C_{n-1}=(2n-1)C_{n-1}+nC_{n-1}=\binom{2n-1}{n}+\binom{2n-2}{n-1}.$$
Thus, to establish \eqref{symvae1}, it suffices to show that its right-hand side gives the cardinality of the set $\mathcal{P}_{(2n-1,1)}\cup\mathcal{P}_{(2n-2,0)}$.  Note that there are $\binom{2n-2}{n-1}$ members of $\mathcal{P}_{(2n-1,1)}$ ending in $u$ and $\binom{2n-2}{n-2}$ ending in $d$.  So to complete the proof of \eqref{symvae1}, one can show
\begin{equation}\label{symvae2}
\binom{2n-1}{n}=1+\sum_{i=1}^{n-1}\left[\binom{2i}{i-1}+\binom{2i-1}{i}\right], \qquad n \geq 2.
\end{equation}

For \eqref{symvae2}, let $\lambda \in \mathcal{P}_{(2n-1,1)}$, with $\lambda\neq (ud)^{n-1}u$.  Consider the largest $i$ such that $\lambda$ passes through either the point  $(2i,2)$ or $(2i-1,-1)$, with $\lambda\neq (ud)^{n-1}u$ implying such an $i \in [n-1]$ exists.  If $\lambda$ passes through $(2i,2)$, then we have the decomposition $\lambda = \alpha d(du)^{n-i-1}$, with $\alpha \in \mathcal{P}_{(2i,2)}$, whereas if $\lambda$ passes through $(2i-1,-1)$, then we have  $\lambda=\beta u(ud)^{n-i-1}u$, with $\beta \in \mathcal{P}_{(2i-1,-1)}$. Thus, for each $i$, there are $\binom{2i}{i-1}+\binom{2i-1}{i}$ possibilities for $\lambda$.  Considering all possible $i \in [n-1]$ implies \eqref{symvae2}, and hence \eqref{symvae1}, as desired. \hfill \qed \medskip

\noindent$\ell$-\emph{valleys:} By an $\ell$-\emph{valley}  within a member of $\mathcal{C}_n$, it is meant a string of the form $ab^\ell(b+1)$ for some $a>b\geq1$ and $\ell \geq1$.  We show that the number of $\ell$-valleys in all the members of $\mathcal{C}_n$ is equal to $\binom{2n-2\ell-1}{n-\ell-3}$ for $1 \leq \ell \leq n-3$, and hence the total number of valleys in $\mathcal{C}_n$ is given by $\sum_{\ell=1}^{n-3}\binom{2n-2\ell-1}{n-\ell-3}$.  Upon adding $\ell-1$ extra copies of the middle letter within a 1-valley so as to obtain an $\ell$-valley, it suffices to prove only the $\ell=1$ case, i.e., that there are $\binom{2n-3}{n-4}$ 1-valleys  in $\mathcal{C}_n$ for all $n \geq 4$.

Note that a 1-valley within $\pi \in \mathcal{C}_n$ corresponds to a string of steps of the form $\tau=ud^ju^2$ in $\imath(\pi)$, where $j \geq 2$.  We thus enumerate the occurrences of $\tau$ in $\mathcal{D}_n$.  To do so, let $\mathcal{J}_n\subseteq \mathcal{D}_{n-2}(u)$ and $\mathcal{K}_n \subseteq \mathcal{D}_{n-2}(d)$ denote the subsets in the respective sets consisting of those members in which the marked $u$ or the marked $d$ is of height at least two.  Let $\alpha=\alpha'{\bf u}\alpha'' \in \mathcal{J}_n$ and $\beta=\beta' ud^i{\bf d}\beta'' \in \mathcal{K}_n$, where the marked step is in bold in each case and $i \geq0$.  Let $\alpha^*=\alpha'{\bf ud^2u^2} \alpha'' \in \mathcal{D}_n(ud^2u^2)$ and $\beta^*=\beta'{\bf ud^{i+3}u^2}\beta'' \in \mathcal{D}_n(ud^{i+3}u^2)$, where the marked string $\rho$ within a member of $\mathcal{D}_n(\rho)$ is in bold here and throughout.  Upon considering all possible $\alpha^*$ and $\beta^*$, one obtains every occurrence in $\mathcal{D}_n$ of $\tau$ for some $j \geq 2$.

Let $\mathcal{J}_n'$ denote the subset of $\mathcal{P}_{(2n-4,2)}$ whose members dip below the $x$-axis and $\mathcal{K}_n'$ the subset of $\mathcal{P}_{(2n-4,-2)}$ whose members dip below the line $y=-3$ at some point.  We may define bijections between $\mathcal{J}_n$ and $\mathcal{J}_n'$ and between $\mathcal{K}_n$ and $\mathcal{K}_n'$ as follows.  Let $\alpha=\alpha'{\bf u}\alpha''$ and $\beta=\beta' {\bf d}\beta''$ denote members of $\mathcal{J}_n$ and $\mathcal{K}_n$, respectively.  Then the mapping $\alpha \mapsto r(\alpha')ur(\alpha'')$ is seen to define a bijection between $\mathcal{J}_n$ and $\mathcal{J}_n'$, upon considering the rightmost step of minimum height within a member of $\mathcal{J}_n'$.  Note that $\alpha \in \mathcal{J}_n$ implies that the marked $u$ is of height at least two, and hence $r(\alpha')$ ends below the $x$-axis (with the endpoint of $r(\alpha')$ being of minimum height).  Similarly, the mapping $\beta \mapsto r(\beta')dr(\beta'')$ is seen to provide a bijection between $\mathcal{K}_n$ and $\mathcal{K}_n'$.

Upon reflecting in the line $y=-1$ all points to the right of and including the first point where a path in $\mathcal{J}_n'$ intersects the line $y=-1$, one can obtain uniquely an arbitrary member of $\mathcal{P}_{(2n-4,-4)}$, by the reflection principle, whence $|\mathcal{J}_n'|=\binom{2n-4}{n-4}$. Similarly, reflection in the line $y=-4$ yields a bijection between members of $\mathcal{K}_n'$ and $\mathcal{P}_{(2n-4,-6)}$, whence $|\mathcal{K}_n'|=\binom{2n-4}{n-5}$.  Combining all of the preceding observations, we then have that the number of 1-valleys in $\mathcal{C}_n$ is given by
$$|\cup_{j=2}^{n-2}\mathcal{D}_n(ud^ju^2)|=|\mathcal{J}_n|+|\mathcal{K}_n|=|\mathcal{J}'_n|+|\mathcal{K}'_n|=\binom{2n-4}{n-4}+\binom{2n-4}{n-5}=\binom{2n-3}{n-4},$$
as desired. \hfill \qed \medskip

\subsection{Peaks}In this subsection, we find combinatorial proofs of the formulas for the totals on $\mathcal{C}_n$ of the parameters tracking the number of symmetric or $\ell$-peaks. ~\medskip

\noindent\emph{Symmetric peaks:} A sequence of letters within a Catalan word of the form $a(a+1)^\ell a$ for some $\ell \geq1$, where $a \geq 1$, is known as a \emph{symmetric} $\ell$-\emph{peak}.  We wish to show that the number of symmetric $\ell$-peaks in $\mathcal{C}_n$ is given by $\binom{2n-2\ell-2}{n-\ell-2}$ for $1 \leq \ell \leq n-2$, and hence there are $\sum_{\ell=1}^{n-2}\binom{2n-2\ell-2}{n-\ell-2}=\sum_{k=0}^{n-3}\binom{2k+2}{k}$ symmetric peaks altogether in $\mathcal{C}_n$ for $n \geq3$, which yields the entry in Table \ref{tab1}.  It suffices again to consider only the $\ell=1$ case, and hence we wish to enumerate all occurrences of $a(a+1)a$ in $\mathcal{C}_n$.

Under $\imath$, an occurrence of $a(a+1)a$ for some $a \geq1$ within a member of $\mathcal{C}_n$ corresponds to a string of steps of the form $u^2d^2u$.  Let $\pi=\pi'{\bf u^2d^2u} \pi'' \in \mathcal{D}_n(u^2d^2u)$, where $\pi',\pi''$ are possibly empty.  The mapping $\pi \mapsto r(\pi')ur(\pi'')$ is seen to define a bijection between $\mathcal{D}_n(u^2d^2u)$ and $\mathcal{P}_{(2n-4,2)}$, which may be reversed by considering the position of the rightmost point of minimum height within a member of $\mathcal{P}_{(2n-4,2)}$.  We thus have that the number of 1-symmetric peaks in $\mathcal{C}_n$ is given by
$$|\mathcal{D}_n(u^2d^2u)|=|\mathcal{P}_{(2n-4,2)}|=\binom{2n-4}{n-3},$$
as desired.  \hfill \qed \medskip

\noindent $\ell$-\emph{peaks:} By an $\ell$-peak within a member of $\mathcal{C}_n$, it is meant a sequence of letters of the form $a(a+1)^\ell b$ for some $a \geq b \geq 1$ and $\ell \geq1$.  We wish to show that the number of $\ell$-peaks in $\mathcal{C}_n$ is equal to $\binom{2n-2\ell-1}{n-\ell-2}$ for $1 \leq \ell \leq n-2$, and hence the total number of peaks is given by $\sum_{\ell=1}^{n-2}\binom{2n-2\ell-1}{n-\ell-2}=\sum_{k=1}^{n-2}\binom{2k+1}{k-1}$. We need only consider the case $\ell=1$ and show that there are $\binom{2n-3}{n-3}$ 1-peaks in $\mathcal{C}_n$ for all $n \geq 3$.

We will describe an occurrence of $u^2d^2$ within a member $\pi \in \mathcal{C}_n$ in which the $d$'s are part of the final run of $d$ steps of $\pi$ as \emph{terminal}, with all other occurrences of $u^2d^2$ in $\pi$ being \emph{non}-\emph{terminal}.  Let $\mathcal{F}_n$ denote the subset of $\mathcal{D}_n(u^2d^2)$ in which the marked $u^2d^2$ is non-terminal.  We have that a $1$-peak within $\rho \in \mathcal{C}_n$ corresponds to a non-terminal occurrence of $u^2d^2$ in $\imath(\rho)$, and hence the total number of 1-peaks in $\mathcal{C}_n$ is given by $|\mathcal{F}_n|$.  Upon inserting $u^2d^2$ at any point along a path in $\mathcal{D}_{n-2}$ (including possibly the start or end points), it is seen that $|\mathcal{D}_n(u^2d^2)|=(2n-3)C_{n-2}=\binom{2n-3}{n-2}$.  Thus, to complete the proof, it suffices to show
\begin{equation}\label{ell-peake1}
|\mathcal{D}_n(u^2d^2) \backslash \mathcal{F}_n|=\binom{2n-3}{n-2}-\binom{2n-3}{n-3}, \qquad n \geq3.
\end{equation}

To do so, let $\pi=\alpha {\bf u^2d^2}d^j \in \mathcal{D}_n(u^2d^2)\backslash \mathcal{F}_n$, where $j \geq0$.  Observe that the final height of the section $\alpha$ within $\pi$ is $j$, and hence the minimum height in $r(\alpha)$ of $-j$ is attained at its endpoint (and possibly other points as well).  Let $\pi^*=u^jdr(\alpha)$.  Note $\pi^* \in \mathcal{P}_{(2n-3,-1)}$ for all $\pi$, with $\pi^*$ not ever going below $y=-1$.  By reflection in the line $y=-2$, members of $\mathcal{P}_{(2n-3,-1)}$ that dip below $y=-1$ are in one-to-one correspondence with members of $\mathcal{P}_{(2n-3,-3)}$, and thus have cardinality $\binom{2n-3}{n-3}$.  Hence, by subtraction, there are $\binom{2n-3}{n-2}-\binom{2n-3}{n-3}$ paths in $\mathcal{P}_{(2n-3,-1)}$ that do not go below $y=-1$.  Since the mapping $\pi \mapsto \pi^*$ is a bijection, the equality \eqref{ell-peake1} follows, which completes the proof.  \hfill \qed \medskip

\subsection{Runs of descents and weak ascents}\label{subs2.3} In this subsection, we treat the formulas for the total number of runs of descents or weak ascents.  ~\medskip

\noindent\emph{Runs of descents:} By a \emph{run of descents} within a sequence $w=w_1\cdots w_n$, it is meant a maximal string $w_iw_{i+1}\cdots w_{i+s}$ for some $i \geq1$ and $s\geq0$ such that $w_i>\cdots>w_{i+s}$.  Note that the number of runs of descents in $\pi \in \mathcal{C}_n$ is one more than the combined number of occurrences of $u^2$ or $udu$ in $\imath(\pi)$.  Thus, the number of runs of descents in $\mathcal{C}_n$ is equal to the cardinality of the set $\mathcal{D}_n\cup\mathcal{D}_n(u^2)\cup\mathcal{D}_n(udu)$.  Let $\alpha=\alpha'{\bf u^2}\alpha'' \in \mathcal{D}_n(u^2)$ and $\beta=\beta'{\bf udu}\beta'' \in \mathcal{D}_n(udu)$.  Then we have that the mappings $\alpha \mapsto r(\alpha')dr(\alpha'')$ and $\beta \mapsto r(\beta')dr(\beta'')$ yield bijections between $\mathcal{D}_n(u^2)$ and $\mathcal{D}_n(udu)$ and the subsets of $\mathcal{P}_{(2n-1,1)}$ and $\mathcal{P}_{(2n-2,0)}$, respectively, consisting of those members that dip below the $x$-axis.  By the reflection principle, these subsets have respective cardinalities $\binom{2n-1}{n-2}$ and $\binom{2n-2}{n-2}$.

Thus, to complete the proof of the formula given in Table \ref{tab1} for runs of descents, we explain combinatorially the identity
\begin{equation}\label{weakasce1}
\binom{2n}{n}=\binom{2n-2}{n-1}+C_n+\binom{2n-1}{n-2}+\binom{2n-2}{n-2}, \qquad n \geq 3.
\end{equation}
Now $\binom{2n-2}{n-1}+\binom{2n-2}{n-2}=\binom{2n-1}{n-1}$, by the recurrence for binomial coefficients, so to establish \eqref{weakasce1}, we need to show $\binom{2n-1}{n-1}=C_n+\binom{2n-1}{n-2}$.  This can be done by noting that there are $C_{n}$ first-quadrant members of $\mathcal{P}_{(2n-1,1)}$, upon appending a terminal $d$ step, and $\binom{2n-1}{n-2}$ members of $\mathcal{P}_{(2n-1,1)}$ that intersect the line $y=-1$, since they are equinumerous with paths in $\mathcal{P}_{(2n-1,-3)}$, by reflection.  This finishes the combinatorial proof of equation \eqref{weakasce1}, as desired. \hfill \qed \medskip

\noindent\emph{Runs of weak ascents:} By a \emph{run of weak ascents} within $w=w_1\cdots w_n$, we mean a maximal string $w_iw_{i+1}\cdots w_{i+s}$ for some $i \geq1$ and $s\geq0$ such that $w_i\leq\cdots\leq w_{i+s}$.  Note that if $w$ is decomposed as $w=w^{(1)}\cdots w^{(r)}$ for some $r \geq1$ such that each $w^{(j)}$ is a run of weak ascents, then $\max(w^{(j)})>\min(w^{(j+1)})$ for $1 \leq j<r$, by maximality, with $w$ containing no other descents.  Thus, the number of runs of weak ascents is one more than the number of descents of $w$ for all $w$.  Therefore, the total number of weak ascents in $\mathcal{C}_n$ is the sum of $C_n$ with the number of descents.

Note that a descent within $\rho \in \mathcal{C}_n$ corresponds to an occurrence of the string $d^2u$ in $\imath(\rho)$.  Let $\pi=\pi'{\bf d^2u}\pi'' \in \mathcal{D}_n(d^2u)$.  Then the section $\pi'$ having final height at least two implies $\pi \mapsto r(\pi')dr(\pi'')$ defines a bijection between $\mathcal{D}_n(d^2u)$ and the subset of $\mathcal{P}_{(2n-2,-2)}$ whose members have minimum height $\leq -3$.  By reflection, such paths are synonymous with arbitrary members of $\mathcal{P}_{(2n-2,-4)}$, and hence have cardinality $\binom{2n-2}{n-3}$.  To show that the number of runs of weak ascents in $\mathcal{C}_n$ is given by $\binom{2n-2}{n-1}$, we thus need to justify combinatorially the equality $\binom{2n-2}{n-1}=\binom{2n-2}{n-3}+C_n$ for $n \geq 3$.

To do so, let $\mathcal{G}_n$ denote the subset of $\mathcal{P}_{(2n-2,0)}$ whose members have minimum height $\leq -2$.  By reflection, there are $\binom{2n-2}{n-3}$ such paths, and hence to complete the proof, we need to show $|\mathcal{P}_{(2n-2,0)}\backslash \mathcal{G}_n|=C_n$.  That is, lattice paths from $(0,0)$ to $(2n-2,0)$ staying at or above the line $y=-1$ are equinumerous with members of $\mathcal{D}_n$.  Let  $\pi \in \mathcal{P}_{(2n-2,0)}\backslash \mathcal{G}_n$ and first suppose $\pi$ does not go below the $x$-axis.  In this case, we let $f(\pi)=u\pi d$.  Otherwise, decompose $\pi$ as $\pi=\pi^{(0)}du\pi^{(1)}\cdots du\pi^{(t)}$ for some $t\geq1$, where $\pi^{(i)}$ for $0 \leq i \leq t$ are all possibly empty Dyck paths.  Let $f(\pi)=u\pi^{(0)}d\cdots u\pi^{(t)}d$, which is seen to be a member of $\mathcal{D}_n$ containing at least two units.  Note that the mapping $f$ may be reversed by considering whether a member of $\mathcal{D}_n$ contains one or more units.  Hence, $f$ provides a bijection between $\mathcal{P}_{(2n-2,0)}\backslash \mathcal{G}_n$ and $\mathcal{D}_n$, which completes the proof.  \hfill \qed \medskip

\noindent \emph{Remarks:} Consider the bargraph representation $b(\pi)$ of $\pi \in \mathcal{C}_n$, viewed as a first-quadrant lattice path from $(0,0)$ to $(n,0)$ using $(0,1)$, $(0,-1)$ and $(1,0)$ steps.  By an \emph{exterior corner} of $b(\pi)$ of type $hu$ ($dh$), it is meant an occurrence of a $(1,0)$  being directly followed by a $(0,1)$ step (a $(0,-1)$ being followed by $(1,0)$, respectively). In \cite{CMR}, it was found that there are $\binom{2n-1}{n-2}$ and $\binom{2n-2}{n-3}$ corners of type $hu$ and $dh$, respectively, in the bargraphs of all the members of $\mathcal{C}_n$.  Note that an $hu$ corner in $b(\pi)$ corresponds to an ascent in $\pi$, and hence to an occurrence of $u^2$ in $\imath(\pi)$, whereas a $dh$ corner corresponds to a descent in $\pi$, and hence to an occurrence of $d^2u$ in $\imath(\pi)$.  In the proofs above for the parameters tracking the runs of descents and the runs of weak ascents, we saw bijectively that $|\mathcal{D}_n(u^2)|=\binom{2n-1}{n-2}$ and $|\mathcal{D}_n(d^2u)|=\binom{2n-2}{n-3}$.  This implies the stated formulas for the total number of corners in $\mathcal{C}_n$ of type $hu$ and $dh$, respectively.

To round out our discussion on runs, we remark that it was shown using generating functions in \cite{BHH} that the runs of ascents and runs of weak descents statistics on $\mathcal{C}_n$ both have the Narayana distribution (e.g., \cite[A001263]{Sloane}). That is, the number of members of $\mathcal{C}_n$ that have exactly $k$ runs of either kind is given by $N(n,k)=\frac{1}{n}\binom{n}{k}\binom{n}{k-1}$ for $1 \leq k \leq n$.  Hence, the total number of runs of either kind is given by $\binom{2n-1}{n}$, which may be afforded a combinatorial argument comparable to those above.

A quick direct argument for the distribution can be given as follows.  Let $\text{asc}(w)$, $\text{des}(w)$ and $\text{lev}(w)$ denote the number of ascents, descents or levels, respectively, in a sequence $w$.  Then the number of runs of ascents in $\pi$ is given by $1+\text{des}(\pi)+\text{lev}(\pi)$ for all $\pi \in \mathcal{C}_n$, since each descent or level of $\pi$ is seen to start a new run of ascents.  This parameter value then translates to
$$1+(\#\,d^2u)+(\#\,udu)=1+(\#\,\text{valleys})=\#\,\text{peaks}$$
in $\imath(\pi)\in \mathcal{D}_n$ for all $\pi$, where a peak or valley within a Dyck path corresponds to an occurrence of $ud$ or $du$, respectively. Since it is well known (with a combinatorial proof) that the peaks distribution on $\mathcal{D}_n$ is Narayana, so too is the distribution for runs of ascents on $\mathcal{C}_n$. Note further that the runs of weak descents parameter on $\mathcal{C}_n$ is the same as $1+\text{asc}$, and hence the equidistribution of the runs of ascents and weak descents follows from the facts $\text{asc}(\pi)+\text{des}(\pi)+\text{lev}(\pi)=n-1$ for all $\pi \in \mathcal{C}_n$ and the Narayana symmetry property $N(n,k)=N(n,n+1-k)$ for $1 \leq k \leq n$.

\section{Combinatorial proof of Catalan identity}

Replacing $k$ with $k+1$ in \eqref{Cniden}, and using the fact $\binom{n}{k}=\frac{n}{k}\binom{n-1}{k-1}$, we show equivalently
\begin{equation}\label{Cniden2}
C_n=\sum_{k=0}^{\lfloor\frac{n-1}{2}\rfloor}\frac{1}{k+1}\binom{n-1}{k}\binom{n-k-1}{k}2^{n-2k-1}, \qquad n \geq 1.
\end{equation}
Let $\mathcal{D}_{n,k}$ for $0 \leq k \leq \lfloor\frac{n-1}{2}\rfloor$ denote the subset of $\mathcal{D}_n$ containing those members in which there are exactly $k$ occurrences of $d^2u$.  Suppose $\pi \in \mathcal{D}_{n,k}$ has exactly $j$ occurrences of the string $udu$. Note $j \leq n-2k-1$ since at least $2k$ $d$'s are involved in internal runs of $d$ of length two or more, with at least one more $d$ in the final run.  Let $\mathcal{D}_{n,k,j}$ denote the subset of $\mathcal{D}_{n,k}$ whose members contain exactly $j$ occurrences of $udu$, and hence $\mathcal{D}_{n,k}=\cup_{j=0}^{n-2k-1}\mathcal{D}_{n,k,j}$ for each $k$.

We first enumerate the members of $\mathcal{D}_{n,k,0}$.  To do so, first suppose,  more generally, that $\pi\in \mathcal{P}_{(2n,0)}$ starts with $u$, ends with $d$ and contains exactly $k+1$ peaks.  In this case, we may write
\begin{equation}\label{piform}
\pi=(\alpha^{(1)}ud\beta^{(1)})(\alpha^{(2)}ud\beta^{(2)})\cdots(\alpha^{(k+1)}ud\beta^{(k+1)}),
\end{equation}
where $\alpha^{(i)}$ and $\beta^{(i)}$ for each $i \in [k+1]$ denote possibly empty runs of $u$ and $d$ steps, respectively.  Let $a_i=|\alpha^{(i)}|$ and $b_i=|\beta^{(i)}|$, and hence $\sum_{i=1}^{k+1}a_i=\sum_{i=1}^{k+1}b_i=n-k-1$.  Then $\pi$ of the form \eqref{piform} is a member of $\mathcal{D}_{n,k,0}$ if and only if $a_i \geq0$ for $i \in [k+1]$, $b_i\geq1$ for $i \in [k]$ and $b_{k+1}\geq0$, with $\sum_{i=1}^{j}a_i\geq \sum_{i=1}^{j}b_i$ for each $j \in [k+1]$.

To determine the permissible vectors of ordered pairs $((a_1,b_1),\ldots,(a_{k+1},b_{k+1}))$ wherein the $a_i$ and $b_i$ satisfy all of the stated conditions above, we invoke an elementary result known as Raney's lemma (e.g., \cite[p.\ 359]{GKP}).  A permutation $\sigma$ of $[m]$ having the form $r(r+1)\cdots m12\cdots(r-1)$ for some $1 \leq r \leq m$ when expressed in the one-line notation is known as a \emph{cylic shift}; that is,
\[
  \sigma(i) =
  \begin{cases}
    i+r-1, & \text{if } 1 \leq i \leq m-r+1, \\
    i-m+r-1, & \text{if } m-r+2 \leq i \leq m.
  \end{cases}
\]
Raney's lemma states that if $(s_1,\ldots,s_m)$ is a sequence of integers whose sum is 1, then there exists exactly one cyclic shift $\sigma$ of $[m]$ such that each of the partial sums $\sum_{i=1}^js_{\sigma(i)}$ for $1 \leq j \leq m$ is positive.

To apply Raney's lemma in this setting, it is more convenient to consider the sequences $u_i=a_{k+2-i}$
and $v_i=b_{k+2-i}$ for $1 \leq i \leq k+1$.  Then $\pi \in \mathcal{D}_{n,k,0}$, expressed as in \eqref{piform}, is uniquely determined by the  vector ${\bf v}$ of ordered pairs  given by
\begin{equation}\label{vform}
{\bf v}=((u_1,v_1),\ldots,(u_{k+1},v_{k+1})),
\end{equation}
wherein $u_i \geq0$ for all $i$, $v_1\geq0$ and $v_i\geq1$ for $i \in [2,k+1]$, with $\sum_{i=1}^jv_i\geq\sum_{i=1}^ju_i$ for each $j \in [k+1]$ and $\sum_{i=1}^{k+1}u_i=\sum_{i=1}^{k+1}v_i=n-k-1$.  To enumerate the permissible vectors ${\bf v}$, we first proceed by inserting $n-k-1$ $u$ and $n-k$ $d$ steps into $k+1$ slots such that each slot receives a non-negative number of $u$'s and a positive number of $d$'s.  Note that this can be implemented in
$$\binom{n-k-1+k}{k}\binom{n-k-(k+1)+k}{k}=\binom{n-1}{k}\binom{n-k-1}{k}$$
ways, by \cite[p.\ 15]{RS}.

Suppose that the $i$-th slot receives $y_i$ $u$'s and $z_i$ $d$'s, where $y_i \geq0$ and $z_i \geq1$ for all $i$.  Let $d_i=z_i-y_i$ for $1 \leq i \leq k+1$ and note $\sum_{i=1}^{k+1}d_i=1$, since one more $d$ was added to the $k+1$ slots than was a $u$.  This allows us to apply Raney's lemma to the sequence $d_i$ and implies that there is a unique cyclic shift $\sigma$ of $[k+1]$ such that $\sum_{i=1}^jd_{\sigma(i)}>0$ for all $1 \leq j \leq k+1$.  Applying this shift to the vector of ordered pairs $(y_i,z_i)$, and subtracting one from the second coordinate of the resulting first pair, yields the vector
$${\bf w}=((y_{\sigma(1)},z_{\sigma(1)}-1),(y_{\sigma(2)},z_{\sigma(2)}),\ldots,(y_{\sigma(k+1)},z_{\sigma(k+1)})).$$
One may verify that ${\bf w}$ is a permissible vector of the form \eqref{vform}.  Thus, there are $$\frac{1}{k+1}\binom{n-1}{k}\binom{n-k-1}{k}$$ permissible vectors ${\bf v}$, and hence the same number of members of $\mathcal{D}_{n,k,0}$, for each $0 \leq k \leq \lfloor \frac{n-1}{2} \rfloor$.

To enumerate the members $\rho \in \mathcal{D}_{n,k,j}$ for $j \geq1$, we start with a precursor $\rho' \in \mathcal{D}_{n-j,k,0}$ and insert $j$ $ud$ strings into the $n-j$ positions directly preceding the $u$ steps of $\rho'$.  Hence, there are $\binom{j+n-j-1}{n-j-1}=\binom{n-1}{j}$ ways in which to insert the strings.  Note that adding the $ud$ strings as described to $\rho'$ does not change the number of occurrences of $d^2u$ nor does it cause the resulting lattice path to dip below the $x$-axis at some point and hence fail to be a Dyck path. Since a member of $\mathcal{D}_{n,k,j}$ is uniquely determined  by its precursor and then independently by the choice of the positions for the added $ud$ strings, we have
$$|\mathcal{D}_{n,k,j}|=|\mathcal{D}_{n-j,k,0}|\binom{n-1}{j}=\frac{1}{k+1}\binom{n-j-1}{k}\binom{n-j-k-1}{k}\binom{n-1}{j}$$
for all $0 \leq k \leq \lfloor \frac{n-1}{2} \rfloor$ and $0 \leq j \leq n-2k-1$.
Therefore, we have
\begin{align*}
|\mathcal{D}_{n,k}|&=\sum_{j=0}^{n-2k-1}\frac{1}{k+1}\binom{n-j-1}{k}\binom{n-j-k-1}{k}\binom{n-1}{j}\\
&=\frac{1}{k+1}\binom{n-1}{k}\binom{n-k-1}{k}2^{n-2k-1},
\end{align*}
where in the last equality we have used the identity
\begin{equation}\label{binomiden}
\sum_{j=0}^{n-2k-1}\binom{n-j-1}{k}\binom{n-j-k-1}{k}\binom{n-1}{j}=\binom{n-1}{k}\binom{n-k-1}{k}2^{n-2k-1}.
\end{equation}

Identity \eqref{binomiden} can be shown by observing
\begin{align*}
\binom{n-1}{j}
\binom{n-j-1}{k}\binom{n-j-k-1}{k}&=\binom{n-1}{k}\binom{n-k-1}{n-j-k-1}\binom{n-j-k-1}{k}\\
&=\binom{n-1}{k}\binom{n-k-1}{k}\binom{n-2k-1}{n-j-2k-1},
\end{align*}
and summing the last product over all $0 \leq j \leq n-2k-1$.  Alternatively, one may provide a quick combinatorial explanation of \eqref{binomiden} as follows.   Note that the right side of \eqref{binomiden} enumerates the set $S$ of ordered triples $(A,B,C)$ consisting of pairwise disjoint subsets $A,B,C$ of $[n-1]$ such that $|A|=|B|=k$, upon choosing the elements for $A$ in $\binom{n-1}{k}$ ways and subsequently selecting the elements for $B$ and $C$.  The left side of \eqref{binomiden} achieves a count of $S$ by considering the number $j$ of elements comprising $C$, in which case there are $\binom{n-j-1}{k}$, $\binom{n-j-k-1}{k}$ and $\binom{n-1}{j}$ possibilities for $A$, $B$ and $C$, respectively.  Summing over all possible $j$ then implies \eqref{binomiden}.

We have thus shown that the generic summand in \eqref{Cniden2} gives the cardinality of $\mathcal{D}_{n,k}$ for each $k$.  Considering all possible $k$ completes the combinatorial proof of \eqref{Cniden2}, as desired.

\section{Semi-perimeter and area totals}

In this section, we give combinatorial explanations of the formulas in Table \ref{tab1} for the total semi-perimeter and area on $\mathcal{C}_n$.

\subsection{Semi-perimeter}
The \emph{bargraph} of $\pi=\pi_1\cdots\pi_n \in \mathcal{C}_n$, denoted by $b(\pi)$, refers to the polyomino in the $(x,y)$-plane containing $n$ adjacent vertical columns of unit width flush with the $x$-axis and starting at the origin such that the $i$-th column has height $\pi_i$ for $1 \leq i \leq n$. Let $\text{semi}(\pi)$ denote the semi-perimeter of $b(\pi)$, which is defined as half the total perimeter of $b(\pi)$ (including its bottom edge along the $x$-axis). \medskip

\noindent\emph{Proof of total semi-perimeter formula:}  We first show
\begin{equation}\label{semie1}
\text{semi}(\pi)=n+1+\text{asc}(\pi), \qquad \pi \in \mathcal{C}_n,
\end{equation}
where $\text{asc}(\pi)$ denotes the number of ascents of $\pi$.

To show \eqref{semie1}, it is useful to view the bargraph $b(\pi)$ as a first-quadrant lattice path, which we will denote by $\pi'$, from $(0,0)$ to $(n,0)$ and containing three types of steps, namely, $x=(0,1)$, $y=(0,-1)$ and $z=(1,0)$, such that no $z$ step occurs at height zero.  Then there are $n$ $z$ steps within $\pi'$, each of which is paired with a segment of unit length along the $x$-axis comprising the bottom boundary of $b(\pi)$, and hence there is a contribution of $n$ towards $\text{semi}(\pi)$.  Further, the first $x$ of $\pi'$, taken together with the last $y$ step, contributes one towards the semi-perimeter in any member of $\mathcal{C}_n$.  All other $x$ steps of $\pi'$ end at height greater than one, with each corresponding to an ascent of $\pi$ since $\pi_{i+1}\leq \pi_i+1$ for all $i$.
We may pair each such $x$ step $\delta$ with the leftmost $y$ step $\delta^*$ to the right of $\delta$ and terminating at the starting height of $\delta$. Each of these $\delta/\delta^*$ pairs contributes one towards $\text{semi}(\pi)$, and taken together, they account for the remaining semi-perimeter of $b(\pi)$.  Combining the prior cases concerning the positions of the steps of $\pi'$ then implies formula \eqref{semie1}.

It was seen in Subsection \ref{subs2.3} that the number of ascents in $\mathcal{C}_n$ is given by $|\mathcal{D}_n(u^2)|=\binom{2n-1}{n-2}$.  Considering all $\pi$ in \eqref{semie1}, we thus have that the sum of the semi-perimeter values of all the members of $\mathcal{C}_n$ is given by
$$C_n+nC_n+|\mathcal{D}_n(u^2)|=C_n+\binom{2n}{n-1}+\binom{2n-1}{n-2}.$$
By the fact $\binom{2m}{m}=2\binom{2m-1}{m-1}$ for all $m$, to complete the proof of the entry for semi-perimeter in Table \ref{tab1}, we thus need to show
\begin{align}
\binom{2n+1}{n}&=\binom{2n-1}{n-1}+C_n+\binom{2n}{n-1}+\binom{2n-1}{n-2}=2\binom{2n}{n-1}+C_n, \quad n \geq 2. \label{semie2}
\end{align}
The second equality in \eqref{semie2} may be realized quickly by arguing that its right-hand side enumerates the members $\pi \in \mathcal{P}_{(2n+1,1)}$. Note that there are $\binom{2n}{n-1}$ possibilities for $\pi$ passing through the point $(2n,2)$, so assume $\pi$ passes through $(2n,0)$.  Then there are $C_n$ possibilities if $\pi$ is first-quadrant and $\binom{2n}{n-1}$ otherwise, since members of $\mathcal{P}_{(2n,0)}$ in the latter case are in one-to-one correspondence with paths in $\mathcal{P}_{(2n,-2)}$, by reflection.  This completes the enumeration of $\mathcal{P}_{(2n+1,1)}$, which implies \eqref{semie2} and finishes the proof.  \hfill \qed

\subsection{Area}By the \emph{area} of $\pi \in \mathcal{C}_n$, it is meant the first-quadrant area subtended by $b(\pi)$, i.e., the sum of all the entries of $\pi$.  \medskip

\noindent\emph{Proof of total area formula:}  Given $\pi \in \mathcal{D}_n$, let $\phi(\pi)$ denote the sum of the heights of all the $u$ steps of $\pi$.  Note that $\text{area}(\pi)=\phi(\imath(\pi))$ for each $\pi \in \mathcal{C}_n$, so we seek a formula for the sum of the $\phi$ values of all the members of $\mathcal{D}_n$. To capture this total using a marked structure, we proceed as follows.  Let us mark any one of the $u$ steps of $\pi \in \mathcal{D}_n$ and suppose that this marked $u$ has height $m$ for some $1 \leq m \leq n$.  We then decompose $\pi$ accordingly as
\begin{equation}\label{areae1}
\pi=\alpha {\bf u}\alpha^{(0)}d\alpha^{(1)}d\cdots \alpha^{(m-1)}d\alpha^{(m)},
\end{equation}
where the marked $u$ is in bold and the $\alpha^{(i)}$ for $0 \leq i \leq m$ are all possibly empty Dyck paths.  Let $\mathcal{D}_n^*$ denote the structure obtained by marking a $u$ step within an arbitrary member of $\mathcal{D}_n$ and choosing any $j \in \{0,1,\ldots,m-1\}$, where $m$ denotes the height of the marked $u$.  Then it is seen that the total $\phi$ value on $\mathcal{D}_n$ is equal to the cardinality of $\mathcal{D}_n^*$, as each $u$ step of height $m$ contributes $m$ towards this cardinality.

To determine $|\mathcal{D}_n^*|$, we define the following operation on $\pi$, expressed as in \eqref{areae1}.  Define $g$ on $\mathcal{D}_n^*$ by putting
\begin{equation}\label{areae2}
g(\pi)=\alpha^{(0)}d\alpha^{(1)}\cdots d\alpha^{(j)}dr(\alpha)dr\left(\alpha^{(j+1)}d\cdots\alpha^{(m-1)}d\alpha^{(m)}\right).
\end{equation}
Note that since the $\alpha^{(i)}$ for all $i$ are Dyck paths, the final height of the initial section $\alpha^{(0)}d\alpha^{(1)}\cdots d\alpha^{(j)}$ is $-j$.  Since the section $\alpha$ of $\pi$ in \eqref{areae1} ends at height $m-1$, we have that the final height achieved after traversing the section $dr(\alpha)d$ of $g(\pi)$ is $-m-j-1$.  Then $r\left(\alpha^{(j+1)}d\cdots\alpha^{(m-1)}d\alpha^{(m)}\right)$ is seen to contribute a net height of $m-j-1$, and hence the final height of $g(\pi)$ is given by $-m-j-1+(m-j-1)=-2j-2$.  Since $0 \leq j \leq m-1$, one obtains lattice paths with endpoints of heights $-2,-4,\ldots,-2m$ for each $\pi \in \mathcal{D}_n^*$ whose marked step is of height $m$.

Let $\mathcal{L}_n$ denote the set of lattice paths with $2n$ steps starting from the origin and having final height negative.  By symmetry, upon excluding members of $\mathcal{P}_{(2n,0)}$, we have
$|\mathcal{L}_n|=\frac{1}{2}\left(4^n-\binom{2n}{n}\right)$. By the preceding observations, we have $g(\pi) \in \mathcal{L}_n$ for all $\pi \in \mathcal{D}_n^*$.  To reverse $g$, suppose that the final and minimum heights achieved by $\lambda \in \mathcal{L}_n$ are $-2j-2$ and $-m-j-1$, respectively, for some $0 \leq j<m\leq n$.  Then consider the leftmost positions $a$ and $b$ at which $\lambda$ achieves a height of $-j-1$ and $-m-j-1$, respectively.  Note that for $g(\pi)$ as in \eqref{areae2}, this occurs with the $d$ steps that directly precede and follow the section $r(\alpha)$.  Once $a$ and $b$ have been specified, it is possible to reconstruct the (unique) pre-image of $\lambda$ with respect to $g$ for any $\lambda$, and hence $g$ is reversible.  Thus $g$ provides a bijection between $\mathcal{D}_n^*$ and $\mathcal{L}_n$, which implies that the total area of $\mathcal{C}_n$ is given by $|\mathcal{L}_n|$, as desired.  \hfill \qed 
\end{document}